\documentclass[letterpaper,11pt,reqno]{amsart}

\usepackage{hyperref}

\usepackage[margin=1in]{geometry}

\theoremstyle{plain}

\newtheorem{theorem}{Theorem}

\theoremstyle{definition}

\theoremstyle{remark}

\newcommand{\BR}{\mathbb{R}}

\newcommand{\BZ}{\mathbb{Z}}

\newcommand{\Leb}{\text{Leb}}

\begin{document}

\title{Brownian bricklayer: a random space-filling curve}
\author{Noah Forman}
\address{Department of Mathematics, University of Washington,
  Seattle, WA 98118, USA}
\email{noah.forman@gmail.com}
\date{\today}
\keywords{Local times, space-filling curve, Brownian motion}
\subjclass[2010]{60J55, 60G17}

\begin{abstract}
 Let $(B(t),\,t\ge0)$ denote the standard, one-dimensional Wiener process and $(\ell(y,t);\,y\in\BR,\,t\ge0)$ its local time at level $y$ up to time $t$. Then $\big((B(t),\,\ell(B(t),t)),\;t\ge0\big)$ is a random path that fills the upper half-plane, covering one unit of area per unit time.
\end{abstract}

\maketitle

The existence of space-filling curves was first proved by Peano, who gave a deterministic, recursive construction, which was later simplified by Hilbert \cite{Mandelbrot77}. In this note we examine a stochastic process, based on Brownian local times, that gives rise to a random space-filling curve. More broadly, this note belongs to the genre of so-called ``pathological'' examples that arise naturally in the study of Brownian motion.  
%
Before presenting the main result, we describe its simpler discrete analogue.


\subsection*{The discrete random bricklayer}

Imagine a worker building a wall by stacking blocks vertically, in side-by-side columns. The worker starts by setting a single block in front of themself. They then flip a coin, stepping to the right if they get heads or to the left if tails before placing the next block. They repeat this process ad infinitum, flipping, moving, and stacking blocks.

We label sites by integers, with the worker's initial position labeled `0'. The worker's left-to-right motion is the \emph{simple random walk on $\BZ$}. The number of blocks stacked at site $j$ after the first $n$ steps is called the \emph{occupation time} of the walk at $j$ up to time $n$.

The simple random walk on $\BZ$ is transitive and recurrent \cite{Durrett}, meaning that it visits every site infinitely many times. Thus, the wall ultimately grows infinitely high at every site. Now suppose the blocks are labeled by the time at which they were placed and we run a string across the front of the infinite wall, connecting consecutively numbered blocks. It is possible that adjacent columns can be at significantly different heights, so the string may run nearly vertical. But this will be uncommon, with the string usually taking smaller vertical steps, since the worker will have visited adjacent sites similar numbers of times \cite[Theorem 1]{CsorReve85}.

\subsection*{The Brownian bricklayer}

For each $n\ge0$, $j\in\BZ$, let $S(n)$ denote the worker's position after $n$ steps and $L(j,n)$ the height of the column at site $j$ at this time. Donsker's theorem states that simple random walk converges in distribution, as a process, in a \emph{scaling limit} to standard, one-dimensional Brownian motion:
\begin{equation}
 \left(n^{-1/2}S\lceil nt\rceil,\,t\ge0\right) \stackrel{d}{\rightarrow} (B(t),\,t\ge0).
\end{equation}
Knight \cite{Knight63} strengthened this result, showing the joint convergence
\begin{equation}
\begin{split}
 &\left(\left( n^{-1/2}S\big(\lceil nt\rceil\big),\,t\ge0 \right),\ \left( n^{-1}L\big(\lceil ny\rceil,\lceil nt\rceil\big);\,t\ge0,\,y\in\BR\right)\right)\\
 &\qquad\stackrel{d}{\rightarrow}\ \big((B(t),\,t\ge0),\ \left(\ell(y,t);\,t\ge0,\,y\in\BR\right)\big),
\end{split}
\end{equation}
where $\ell(y,t)$ is the (occupation density) \emph{local time} of Brownian motion at level $y$, up to time $t$. This is defined as
\begin{equation}\label{eq:LT_def}
 \ell(y,t) := \lim_{\epsilon\to 0}(2\epsilon)^{-1}\Leb\{u\in [0,t]\colon |B(u)-y| < \epsilon\} \qquad \text{for }y\in\BR,\,t\ge0.
\end{equation}
Informally, the local time $\ell(y,t)$ quantifies the amount of time spent by the Brownian path at (or more accurately, near) level $y$, prior to time $t$. See \cite[Chapter 6]{MortersPeres} for more background. 

We define the \emph{Brownian bricklayer} to be the process
\begin{equation}
 K(t) := \big(B(t),\;\ell(B(t),t)\big),\qquad t\ge 0.
\end{equation}
In analogy with our discrete-time description, the first coordinate of this process describes the current location of a bricklayer as they move up and down the line, and the second describes the height of the wall in front of the bricklayer's current location. Extending the analogy, the wall itself at time $t$ is described by the \emph{local time profile} $(\ell(y,t),\,y\in\BR)$.

\begin{theorem}\label{thm:main}
 The process $\big(K(t),\,t\ge0\big)$ is almost surely (a.s.) path continuous, and its path a.s.\ fills the upper half plane. I.e.\ $\left\{K(t) \colon t\ge 0\right\} = \BR\times [0,\infty)$ almost surely.
\end{theorem}

To prove this result, we appeal to two well-known properties of Brownian local times.
\begin{enumerate}
 \item Trotter's theorem \cite[Theorem 6.19]{MortersPeres}: the random field $(\ell(y,t);\, y\in\BR,t\ge0)$ admits an a.s.\ continuous version. I.e.\ there is an a.s.\ event on which the limit $\ell(y,t)$ is defined and continuous at all points $(y,t)\in \BR\times [0,\infty)$.
 \item It is a.s.\ the case that, for every $y\in\BR$, the process $(\ell(y,t),\,t\ge0)$ increases without bound. This can be deduced from Ray's Theorem \cite[Theorem 6.38]{MortersPeres} and the strong Markov property.
\end{enumerate}

Let $\Omega'$ denote an a.s.\ event on which both of the above properties hold.

\begin{proof}[Proof of Theorem \ref{thm:main}]
 The continuity of $(K(t),\,t\ge0)$ follows from Trotter's theorem and the continuity of Brownian motion itself. It remains to show that the path is a.s.\ space-filling. To that end, fix $(y,s)\in\BR\times [0,\infty)$. It suffices to show that for every outcome $\omega\in \Omega'$ there is some $T = T(\omega,y,s)$ for which $\big(B(T),\ell(B(T),T)\big) = (y,s)$.
 
 If $s=0$ then, by the continuity of $\ell$, taking $T = \inf\{t\ge 0\colon B(t) = y\}$ gives the desired result. Now, assume $s>0$. Consider the set $\{t\ge 0\colon \ell(y,t) < s\}$. This set is non-empty and, by property (ii) above, it is bounded, so it has a supremum $T$. By the continuity of $\ell$ we get $\ell(y,T) = s$, as desired. Moreover, for all $t<T$ we have $\ell(y,t) < \ell(y,T)$. By the definition in \eqref{eq:LT_def}, this implies that for every $\epsilon>0$ and every $t<T$, there is some $t'\in (t,T)$ for which $|B(t')-y| < \epsilon$. In particular, by the continuity of $B$ we get $B(T) = y$, also as desired.
\end{proof}

\subsection*{Properties of the Brownian bricklayer}

The following results are consequences of well known properties of Brownian local times and excursions. See \cite{BertoinLevy,GreePitm80,MortersPeres} for general references on these topics. 
\smallskip

(a) \emph{The area of the brick wall grows determinstically, at unit rate}. For each $y$, the local time $(\ell(y,u),\,u\ge 0)$ increases only at times $u$ when $B(u) = y$. Thus, at time $t$, at each position $y$, the wall reaches height $\ell(y,t)$. I.e.
 $$\big\{K(u),\, u\in [0,t]\big\} = \big\{ (y,s)\colon y\in \{B(u),\,u\in [0,t]\},\ s\in [0,\ell(y,t)] \big\}.$$
 %
 The area of the wall at time $t$ is thus $\int_{-\infty}^{\infty} \ell(y,t)dy = t$, by the definition in \eqref{eq:LT_def}.
 
 For $c\neq 0$ and $d>0$, the process $\big(\big(cB(t),d\ell(B(t),t)\big),\ t\ge0\big)$ is a space-filling curve covering an area that grows deterministically, at rate $|c|d$.\smallskip

(b) \emph{The process $(K(t),\,t\ge0)$ is non-Markovian}. Its future evolution depends on the local time profile, $(\ell(y,t),\,y\in\BR)$, which we think of as the state of the growing wall. However, the augmented process 
$\big(\big(K(t),\big(\ell(y,t),\,y\in\BR\big)\big),\ t\ge0\big)$ 
is a strong Markov process.\smallskip

(c) Barlow \cite{Barlow82} proved that $\big(\ell(B(t),t),\,t\ge0\big)$ is \emph{not a semimartingale} by showing that it fails to be H\"older-continuous with index $\frac14$.\smallskip

(d) The law of $K(t)$, at fixed times $t$, is 
\begin{equation}\label{eq:one_dim_dist}
 \Pr\left\{B(t)\in dy,\, \ell(B(t),t)\in ds\right\} = \frac{|y|+s}{\sqrt{8\pi t^3}}\exp\left(-\frac{(|y|+s)^2}{2t}\right)dyds.
\end{equation}
We deduce this by a sequence of transformation identities. First, the identity $(B(u),\,u\in [0,t]) \stackrel{d}{=} (B(t)-B(t-u),\,u\in [0,t])$ implies that
\begin{equation}
 \left(B(t),\ell(B(t),t)\right) \stackrel{d}{=} \left(B(t),\ell(0,t)\right).
\end{equation}
L\'evy's theorem \cite[Theorem 7.38]{MortersPeres} states that
\begin{equation}
 \left(|B(t)|,\ell(0,t)\right) \stackrel{d}{=} (S(t)\!-\!B(t),\,S(t)), \quad \text{where} \quad S(t) := \max_{u\in [0,t]}B(u).
\end{equation}
The identity $(B(u),\,u\in [0,t]) \stackrel{d}{=} (-B(u),\,u\in [0,t])$ implies that the sign of $B(t)$ is independent of $\left(|B(t)|,\ell(0,t)\right)$. Thus,
\begin{equation}
 \left(B(t),\ell(0,t)\right) \stackrel{d}{=} ((S(t)-B(t))\cdot I,S(t)),
\end{equation}
where $I = \pm1$ with probability $\frac12$ each, independent of the Brownian motion. The reflection principle (see \cite{MortersPeres}) implies that, for $x\in\BR$ and $s>0$ with $x<s$,
\begin{equation}
 \Pr\left\{B(t) \in dx,\, S(t)\ge s\right\} = \Pr\{2s - B(t) \in dx\} = \frac{1}{\sqrt{2\pi t}}\exp\left(-\frac{(2s-x)^2}{2t}\right).
\end{equation}
Putting these pieces together gives \eqref{eq:one_dim_dist}.\smallskip

(e) For any fixed $y$, the set 
of points on the line $\{y\}\times [0,\infty)$ that the process $K$ visits multiple times is a.s.\ countable, and no point on this line gets visited more than twice. In particular, $\{t\colon B(t)=y\}$ is known to be homeomorphic to the Cantor middle-third set, with the removed middle thirds corresponding to excursions of $B$ away from $y$. If $(G,D)$ denotes the time interval of one such excursion, then $B(G) = B(D) = y$ and $\ell(y,G) = \ell(y,D)$. Thus $K(t)$ revisits a previous value in $\{y\}\times [0,\infty)$ at the end of each excursion. But there are a.s.\ no two excursions about level $y$ arising at the same local time, so these are a.s.\ the only repetitions.

\subsection*{Related literature}

As noted above, Barlow \cite{Barlow82} proved that $\big(\ell(B(t),t),\,t\ge0\big)$ is not a semimartingale. Subsequent papers \cite{Aldous86,Bertoin89} have encountered this process as well. It is closely related to $\big(\sup_y\ell(y,t),\,t\ge0)$, since 
$\sup_y\ell(y,t) = \sup_{u\in [0,t]}\ell(B(u),u)$. 
See \cite{Revesz}, for example, for discussion of this supremum.

The local time profile $(\ell(y,T),\,y\in\BR)$ at a deterministic or random time $T$ has been studied extensively, most notably with the Ray-Knight theorems \cite[Chapter 6]{MortersPeres}. 
These theorems state that, at certain natural random times, the local time as a process in level behaves like a certain continuous state branching process, which is a natural model for random population growth and decay. Jeulin's theorem \cite{Jeulin85} gives an alternative descrition of the local time profile of a Brownian excursion, as a certain time-change of another Brownian excursion. Generalizations of the Ray-Knight and Jeulin descriptions to fixed times appear in \cite{Leuridan98} and \cite{AssaFormPitm15,LupuPitmTang13}, respectively.

The name ``Brownian bricklayer'' is an homage to Warren and Yor's ``Brownian burglar'' \cite{WarrYor98}: Brownian motion conditioned on its local time profile. 
Aldous \cite{Aldous98} conducted a similar study. The local time profile, as a process evolving in time, arises in the study of self-repelling processes, such as true self-repelling motion (TSRM) \cite{TothWern98}. See \cite{Pemantle07} for a survey of foundational results in this area.

On the discrete side, many authors have studied the occupation times, or discrete local times, of simple random walks. Such counts were the basis for Knight's proof of the Ray-Knight theorems \cite{Knight63}. See also subsequent work by Cs\"org\H{o} and R\'ev\'esz, such as \cite{CsorReve84,Revesz81}.

Schramm-Loewner evolution with parameter $\kappa > 8$ (SLE$(\kappa)$) is another random path that fills the upper half-plane \cite{RohdSchr05}, but it does not equal the Brownian bricklayer. For example, the bricklayer fills each vertical column from bottom to top -- it cannot visit the point $(0,1)$ before visiting $(0,1/2)$, whereas SLE can.

\subsection*{Acknowledgements}
The author thanks Soumik Pal, Jim Pitman, and Matthias Winkel for their helpful comments.

\bibliographystyle{plain}
\bibliography{SpaceFilling}
\end{document}